\input amstex
\documentstyle{amsppt}
\input epsf
\magnification 1200
\vcorrection{-1cm}
\NoBlackBoxes

\def\R{\Bbb R}
\def\C{\Bbb C}

\def\GF{\Bbb F}
\def\tr{\operatorname{tr}}
\def\diag{\operatorname{diag}}
\def\smat{\left(\smallmatrix}
\def\esmat{\endsmallmatrix\right)}
\def\eps{\varepsilon}

 \let\Cl=\cl
\def\co{{\text{c}}}

\def\sectClasses {2}
\def\sectMainRes {3}
\def\sectTr      {4}
\def\sectDoubleProd {5}
\def\sectTripleProd {6}
\def\sectQuadrProd  {7}

\def\propClass {2.1}

\def\thMain   {3.1}
\def\remGplus {3.2}
\def\corPSL   {3.3}
\def\corCN    {3.4}

\def\lemSU       {4.1}
\def\lemEV       {4.2}
\def\lemAB       {4.3}
\def\lemTwoThree {4.4}
\def\lemTwoTwo   {4.5}
\def\lemCFour    {4.6}

\def\lemABG        {5.1}
\def\lemTwoThreeAB {5.2}
\def\lemTwoTwoAB   {5.3}
\def\lemCFourAB    {5.4}

\def\eqIsom  {1}
\def\eqABone {2}

\def\tabClasses  {1}
\def\tabProd     {2}
\def\tabProdABG  {3} 

\def\figABG  {1}

\def\refAW      {1}                       
\def\refAH      {2}                       
\def\refACC     {3}                       
\def\refBelkale {4}
\def\refB       {5}                       
\def\refFW      {6}                       
\def\refGP      {7}                       
\def\refGordeev {8}
\def\refKarni   {9} \let\refK=\refKarni
\def\refLSh     {10}
\def\refJKTR    {11}
\def\refAFST    {12}
\def\refP       {13}                     
\def\refPW      {14}                     
\def\refS       {15}                     

\topmatter
\title         Products of conjugacy classes in $SL_2(\R)$
\endtitle
\author        S.~Yu.~Orevkov
\endauthor
\abstract      We compute the product of any $n$-tuple
               of conjugacy classes in $SL_2(\R)$.
\endabstract
\endtopmatter

\document

\head 1. Introduction
\endhead
In this paper we compute the product of any $n$-tuple
of conjugacy classes in the group $SL_2(\R)$ (see Theorem \thMain\ in \S\sectMainRes).
The computation is straightforward, the main
difficulty being to find a suitable notation for the answer 
to be readable.

All products of conjugacy classes are computed in
[\refK]  for simple finite groups of order less than
million and for sporadic simple finite groups. In 
[\refAFST] the same is done 
for finite unitary groups $GU_3(\GF_{q})$ and $SU_3(\GF_{q})$
as well as for finite linear
groups $GL_3(\GF_{q})$ and $SL_3(\GF_q)$.
Similar questions were studied by several authors, see
[\refAH, \refB, \refGordeev, \refLSh, \refS] and references therein.

My special interest in the computation of class products in any kind of
linear or unitary groups is motivated by possible applications to
plane real or complex algebraic curves, see [\refACC, \refJKTR].
Perhaps, the most interesting and
non-trivial case when the products of conjugacy classes are completely computed, is
the case of the unitary groups $SU(n)$, see [\refAW, \refBelkale].
It seems that Belkale's approach [\refBelkale] could be extended (at least
partially) to pseudo-unitary groups $SU(p,q)$ using the techniques developed in
[\refGP]. We are going to do this in a subsequent paper. Some products
of conjugacy classes in $PU(n,1)$ (especially for $n=2$) are computed in
[\refFW, \refP, \refPW].

Note that $SU(1,1)$ is isomorphic to $SL_2(\R)$. Indeed,
$$
  \Phi:SL_2(\R)\to SU(1,1),\qquad
   A \mapsto\Phi(A)= P^{-1}AP, \qquad
   P=\left(\matrix 1 & i\\i & 1\endmatrix\right)\,,        \eqno(\eqIsom)
$$
is an isomorphism; recall that $SU(1,1)=\{A\in SL_2(\C)\mid A^*JA=J\}$
where $J=\diag(1,-1)$. So, the main motivation for the computation of
class products in $SL_2(\R)$ was to get an idea what the answer for $SU(p,q)$
could look like.

\head \sectClasses. Conjugacy classes
\endhead
Let $G=SL_2(\R)$.
The conjugacy classes in $G$ are given in Table \tabClasses.
This fact can be easily derived, e.~g., from [\refB, \S2].

\midinsert
\centerline{ \vbox{\offinterlineskip
\def\h {height6pt&\omit&&\omit&&\omit&&\omit&\cr}
\def\o{\omit} \def\e{\varepsilon} \def\d{\delta} \def\a{\alpha} \def\la{\lambda}
\def\s {\left(\smallmatrix} \def\sm{\setminus}
\def\es{\endsmallmatrix\right)} \def\and{\quad\&\quad}
\hrule
\halign{&\vrule#&\strut\;\hfil#\hfil\,\cr
\h
&  \lower-4pt\hbox{\;class}
&& \lower-4pt\hbox{\;paramameters}
&& \lower-4pt\hbox{\;representative}
&& \hbox{\vbox{\pagewidth{50mm} nessessary and sufficient\par condition on
     $A=\s a&b\\c&d\es$ }}&\cr
\h
 \noalign{\hrule}
\h
& $\Cl_1^\e$      && $\{-1,1\}$      && $\e I=\s\e&0\\ 0&\e\es$   && $A=\e I$ &\cr\h
& $\Cl_2^{\e,\d}$ && $\{-1,1\}^2$    && $\s\e&0\\\d&\e\es$  && $\tr A=2\e\and\d(c-b)>0$ &\cr\h
& $\Cl_3^\a$ && $]0,2\pi[\,\sm\{\pi\}$ && $\s\cos\a&-\sin\a\\
                                       \sin\a& \cos\a\es$ && $\tr A=2\cos\a\and c\sin\a>0$ &\cr\h
& $\Cl_4^\la$     && $\R\sm[-1,1]$         && $\s\la&0\\
                                                0&\la^{-1}\es$ && $\tr A=\la+\la^{-1}$ &\cr\h
\noalign{\hrule}
}} }
\botcaption{ Table \tabClasses } Conjugacy classes in $SL_2(\R)$.
Note that for $\Cl_2^{\eps,\delta}$ (resp.  $\Cl_3^\alpha$) we have
$\delta(c-b)>0\,\Leftrightarrow\,(\delta c>0\text{ or }\delta b<0)$
(resp.
$c\sin\alpha>0\,\Leftrightarrow\,b\sin\alpha<0$)
\endcaption
\endinsert

A more geometric characterization of the conjugacy classes can be given
as follows. For $\vec x=(x_1,x_2)$ and $\vec y=(y_1,y_2)$, we denote
$\vec x\wedge\vec y=x_1y_2-x_2y_1$.

\proclaim{ Proposition \propClass } Let $A\in G\setminus\{I,-I\}$ and $0<\alpha<\pi$. Then:

(a). $A\in\Cl_2^{\eps,\delta}$, $\eps,\delta=\pm1$,
if and only if $\tr A=2\eps$ and $\delta\vec x\wedge A\vec x\ge 0$
for any $\vec x\in\R^2$;

(b). $A\in\Cl_2^\alpha$ if and only if $\tr A=2\cos\alpha$ and  $\vec x\wedge A\vec x>0$
for any $\vec x\in\R^2$. \qed
\endproclaim

We denote the union of all $\Cl_i^{\dots}$ by $\Cl_i$ ($i=1,2,3,4$).
For $X\subset G$, we denote the complement $G\setminus X$ by $X^\co$.
We set also
$$
   \Cl_2^{\eps,*}=\Cl_2^{\eps,+}\cup\Cl_2^{\eps,-},\qquad
   \Cl_4^+=\bigcup_{\lambda>1} \Cl_4^\lambda,\qquad
   \Cl_4^-=\bigcup_{\lambda<-1} \Cl_4^\lambda,
$$
$$
   \overline{\Cl_4^+}=\Cl_4^+\cup\{ I\}\cup \Cl_2^{+*},
\qquad
   \overline{\Cl_4^-}=\Cl_4^-\cup\{-I\}\cup \Cl_2^{-*}.
$$
If $j$ is an interval (open or closed from either end)
contained in $]0,2\pi[\setminus\{\pi\}$, then we set
$\Cl_3^j = \bigcup_{\alpha\in j} \Cl_3^\alpha$.
If the left end of $j$ is ``$]0\dots$'' or
``$]\pi\dots$'', then the bracket ``$]\dots$''
can be replaced by ``$[\dots$'' or ``$\langle[\dots$''. Symmetrically,
if the right end of $j$ is ``$\dots\pi[$'' or
``$\dots2\pi[$'', then the bracket ``$\dots[$''
can be replaced by ``$\dots]$'' or ``$\dots]\rangle$'',
and this means the following:
$$
\xalignat2
   &\Cl_3^{[0,\dots} = \Cl_2^{++}\cup \Cl_3^{]0,\dots},   &
   &\Cl_3^{\dots,\pi]}=\Cl_3^{\dots,\pi[}\cup \Cl_2^{-+}, \\
   &\Cl_3^{[\pi,\dots}=\Cl_2^{--}\cup \Cl_3^{]\pi,\dots}, &
   &\Cl_3^{\dots,2\pi]}=\Cl_3^{\dots,2\pi[}\cup \Cl_2^{+-},
\\
   &\Cl_3^{\langle[0,\dots} = \Cl_4^{+}\cup \Cl_3^{[0,\dots},   &
   &\Cl_3^{\dots,\pi]\rangle}=\Cl_3^{\dots,\pi]}\cup \Cl_4^{-}, \\
   &\Cl_3^{\langle[\pi,\dots}=\Cl_4^{-}\cup \Cl_3^{[\pi,\dots}, &
   &\Cl_3^{\dots,2\pi]\rangle}=\Cl_3^{\dots,2\pi]}\cup \Cl_4^{+},
\endxalignat
$$
for example,
$\Cl_3^{\langle[0,\pi]}=\Cl_4^+\cup \Cl_2^{++}\cup
   \Big(\bigcup_{0<\alpha<\pi} \Cl_3^\alpha\Big)\cup \Cl_2^{-+}$.
Let
$$
     G^+ = \overline{\Cl_4^+}\cup \Cl_3^{]0,\pi[}.
$$

\head \sectMainRes. Statement of the main result
\endhead

\proclaim{ Theorem \thMain }
(a). For conjugacy classes contained in $G^+$ (see Remark \remGplus), their
double and triple products are as shown in Tables \tabProd\ and \tabProdABG;
in Figure \figABG\ we represent all triples of classes from $\Cl_3$
whose product contains $I$.

\smallskip
(b). Let $0<\alpha,\beta,\gamma,\delta<\pi$. Then:
$$
\xalignat2
  &\Cl_2^{++}\Cl_2^{++}\Cl_2^{+-}\Cl_2^{+-}
   = \Cl_2^{++}\Cl_2^{++}\Cl_2^{+-}\Cl_3^\alpha=\{-I\}^\co,\\
  &\Cl_2^{++}\Cl_2^{+-}\Cl_2^{+-}\Cl_3^\alpha
   = \Cl_2^{++}\Cl_2^{+-}\Cl_3^\alpha\Cl_3^\beta=\{I\}^\co,\\
  &\Cl_2^{++}\Cl_2^{++}\Cl_3^\alpha\Cl_3^\beta = \{-I\}^\co &&\text{if $\alpha+\beta\ge\pi$},\\
  &\Cl_2^{+-}\Cl_2^{+-}\Cl_3^\alpha\Cl_3^\beta = \{I\}^\co &&\text{if $\alpha+\beta\le\pi$},\\
  &\Cl_3^\alpha\Cl_3^\beta\Cl_3^\gamma\Cl_3^\delta = \{-I\}^\co &&\text{if
            $\pi<\alpha+\beta+\gamma+\delta<3\pi$},
\endxalignat
$$
the products of the form $\Cl_3^\alpha\Cl_3^\beta\Cl_3^\gamma\Cl_2^{+\pm}$ are as in
Table \tabProdABG, and for any other four non-scalar conjugacy classes
contained in $G^+$, their product is the whole $G$.

\smallskip
(c). The product of any five non-scalar conjugacy classes is the whole $G$.
\endproclaim

\proclaim{ Remark \remGplus }{\rm
We see in Table~\tabClasses\ that, for any conjugacy class $\Cl$. either $\Cl\subset G^+$
or $-\Cl\subset G^+$.
Thus any product of conjugacy classes can be immediately recovered if one knows
all products of classes contained in $G^+$.}
\endproclaim

Let $\tilde G=PSL_2(\R)$. We denote the image of $\Cl_i^{\dots}$ in $\tilde G$
by $\tilde\Cl_i^{\dots}$.


\proclaim{ Corollary \corPSL } (a).
We have $\tilde\Cl_4^\lambda\tilde\Cl_4^\lambda=\tilde G$ and
$$
   \tilde\Cl_2^{++}\tilde\Cl_2^{++} = \tilde\Cl_2^{+-}\tilde\Cl_2^{+-}
    = \tilde\Cl_4^\lambda \tilde\Cl
    = \tilde G\setminus\{\tilde I\},
\qquad\text{where $\tilde\Cl\ne\tilde\Cl_4^\lambda,\tilde I$}
$$
(see Table \tabProd\ for the other double products of classes in $\tilde G$).

\smallskip
(b). Let $0<\alpha,\beta,\gamma<\pi$. The following triple products are equal to
$\tilde G\setminus\{\tilde I\}$:
$$
\split
    &\tilde\Cl_2^{++}\tilde\Cl_2^{+-}\tilde\Cl_3^\alpha,\quad
    \tilde\Cl_2^{+-}\tilde\Cl_3^\alpha\tilde\Cl_3^\beta \;
    \text{(if $\alpha+\beta<\pi$)},\quad
    \tilde\Cl_2^{+\pm}\tilde\Cl_3^\alpha\tilde\Cl_3^\beta \;
    \text{(if $\alpha+\beta=\pi$)},
\\
    &\tilde\Cl_2^{++}\tilde\Cl_3^\alpha\tilde\Cl_3^\beta \;
    \text{(if $\alpha+\beta>\pi$)},\quad
    \tilde\Cl_3^\alpha\tilde\Cl_3^\beta\tilde\Cl_3^\gamma \;
    \text{(if $\pi<\alpha+\beta+\gamma<2\pi$)}.
\endsplit
$$
The product of any other three non-trivial conjugacy classes is the whole $\tilde G$.

\smallskip
(c). The product of any four non-trivial conjugacy classes is the whole $\tilde G$.
\endproclaim

\proclaim{ Corollary \corCN }
$\operatorname{cn}(\tilde G)=
 \operatorname{ecn}(\tilde G)=4$ (in the notation of [\refK]).
\endproclaim

\midinsert
\centerline{ \vbox{\offinterlineskip
\def\h {height2pt&\omit&&\omit&&\omit&&\omit&&\omit&&\omit&\cr}
\def\hh{height4pt&\omit&&\omit&&\omit&&\omit&&\omit&&\omit&\cr}
\def\o{\omit} \def\ra{\rangle} \def\la{\langle}
\def\PP  {$\Cl_3^{[0,\pi]\ra}$}
\def\PPP {$\{I\}^c$}
\def\PPM {$\Cl_3^{[0,\pi]\ra}\cup\overline{\Cl_4^+}$}
\def\PPA {$\big(\{I\}\cup\Cl_3^{[0,\alpha]}\big)^\co$}

\def\PPL {$\{I\}^\co$}
\def\PM  {$\overline{\Cl_4^{+}}$}
\def\PMM {$\Cl_3^{\la[\pi,2\pi]}\cup\overline{\Cl_4^+}$}
\def\PMA {$\Cl_3^{\la[0,\pi]\ra}$}

\def\PML {$\{-I\}^\co$}
\def\MM  {$\Cl_3^{\langle[\pi,2\pi]}$}
\def\MMM {$\{I\}^\co$}
\def\MMA {$\big(\{-I\}\cup\Cl_3^{[\alpha,\pi]}\big)^\co$}

\def\MML {$\{I\}^\co$}
\def\PA  {$\Cl_3^{\,]\alpha,\pi]\ra}$}
\def\MA  {$\Cl_3^{\langle[0,\alpha[}$}
\def\MMA {$\big(\{-I\}\cup\Cl_3^{[\alpha,\pi]}\big)^\co$}
\def\ApP {$\big(\{I\}\cup\Cl_3^{[0,\alpha+\beta]}\big)^\co$}
\def\ApM {$\Cl_3^{\la[0,\pi]\ra}$}
\def\AzP {$\Cl_3^{\la[\pi,2\pi]\ra}$} 
\def\AzM {$\Cl_3^{\la[0,\pi]\ra}$}
\def\AnP {$\Cl_3^{\la[\pi,2\pi]\ra}$}
\def\AnM {$\big(\{I\}\cup\Cl_3^{[\alpha+\beta,\pi]}\big)^\co$}
\def\PAL {$\{I\}^\co$}
\def\MAL {$\{-I\}^\co$}

 \let\AGL\ABL

\def\PL  {$\Cl_3^{\la[0,\pi]\ra}$}
\def\ML  {$\Cl_3^{\langle[\pi,2\pi]\rangle}$}
\def\Ap  {$\Cl_3^{[\alpha+\beta,\pi]\rangle}$}
\def\Az  {$\{-I\}\cup\Cl_4^-$}
\def\An  {$\Cl_3^{\la[\pi,\alpha+\beta]}$}
\def\AL  {$\Cl_3^{\la[0,\pi]\ra}$}
\def\LL  {$\{-I\}^\co$}
\def\Lmu {$\{I,-I\}^\co$}
\def\SW {$\swarrow$}

\let\tabl\tabProdABG
\hrule
\halign{&\vrule#&\strut\,\hfil#\hfil\,\cr
\h
&\omit && $I$ && $\Cl_2^{++}$ && $\Cl_2^{+-}$ && $\Cl_3^\gamma$ && $\Cl_4^\nu$ &\cr
\h
 \noalign{\hrule}
\h
& $\Cl_2^{++}\Cl_2^{++}$      && \PP && \PPP && \PPM && \SW  && \SW  &\cr\hh
& $\Cl_2^{++}\Cl_2^{+-}$      && \PM && \PPM && \PMM && \SW  && \SW  &\cr\hh
& $\Cl_2^{+-}\Cl_2^{+-}$      && \MM && \PMM && \MMM && \SW  && \SW  &\cr\hh
& $\Cl_2^{++}\Cl_3^{\alpha}$  && \PA && \PPA && \PMA && \SW  && \SW  &\cr\hh
& $\Cl_2^{+-}\Cl_3^{\alpha}$  && \MA && \PMA && \MMA && \SW  && \SW  &\cr\hh
& $\Cl_2^{++}\Cl_4^{\lambda}$ && \PL && \PPL && \PML && \SW  && \SW  &\cr\hh
& $\Cl_2^{+-}\Cl_4^{\lambda}$ && \ML && \PML && \MML && \SW  && \SW  &\cr\hh
& $\Cl_3^\alpha\Cl_3^\beta$,
  $\alpha+\beta<\pi$      && \Ap && \ApP && \ApM && \lower8pt\hbox{See}    && \SW  &\cr
& $\Cl_3^\alpha\Cl_3^\beta$,
  $\alpha+\beta=\pi$      && \Az && \AzP && \AzM && \lower8pt\hbox{Tbl.~\tabl} && \SW  &\cr
& $\Cl_3^\alpha\Cl_3^\beta$,
  $\alpha+\beta>\pi$      && \An && \AnP && \AnM && \o   && \SW  &\cr\hh
& $\Cl_3^\alpha\Cl_4^\lambda$&& \AL && \PAL && \MAL && \AGL && $G$  &\cr\hh
&$\Cl_4^\lambda\Cl_4^\lambda$&& \LL && $G$  && $G$  && $G$  && $G$  &\cr\hh
&$\Cl_4^\lambda\Cl_4^\mu$, 
 $\lambda\ne\mu$         && \Lmu && $G$  && $G$  && $G$  && $G$  &\cr\h
\noalign{\hrule}
}}
 }
\botcaption{ Table \tabProd }
Double and triple products of conjugacy classes in $SL_2(\R)$
(``$\swarrow$'' means ``see some other cell(s) of this table"). The range of the parameters:
$0<\alpha,\beta,\gamma<\pi$ and $\lambda,\mu,\nu>1$.
\endcaption
\endinsert

\midinsert
\centerline{ \vbox{\offinterlineskip
\def\h {height6pt&\omit&&\omit&&\omit&&\omit&\cr}
\def\o{\omit} \def\e{\varepsilon} \def\d{\delta} \def\a{\alpha} \def\la{\lambda}
\def\s {\left(\smallmatrix} \def\sm{\setminus}
\def\es{\endsmallmatrix\right)} \def\and{\quad\&\quad}
\hrule
\halign{&\vrule#&\strut\;\hfil#\hfil\,\cr
\h
&  \;Condition on $\alpha,\beta,\gamma$
&& \;$\Cl_3^\alpha \Cl_3^\beta \Cl_3^\gamma$
&& \;$\Cl_3^\alpha \Cl_3^\beta \Cl_3^\gamma \Cl_2^{++}$
&& \;$\Cl_3^\alpha \Cl_3^\beta \Cl_3^\gamma \Cl_2^{+-}$
&\cr
\h
 \noalign{\hrule}
\h
& $\alpha+\beta+\gamma<\pi$  && $\big(\{I\}\cup \Cl_3^{[0,\alpha+\beta+\gamma[}\big)^\co$
  && $G$ && $G\setminus\{I\}$
  &\cr\h
& $\alpha+\beta+\gamma=\pi$  && $\{-I\}\cup \Cl_3^{\langle[\pi,2\pi]\rangle}$
  && $G\setminus\{-I\}$ && $G\setminus\{I\}$
  &\cr\h
& $\pi<\alpha+\beta+\gamma<2\pi$ && $\Cl_3^{\langle[\pi,2\pi]\rangle}$
  && $G\setminus\{-I\}$ && $G\setminus\{I\}$
  &\cr\h
& $\alpha+\beta+\gamma=2\pi$ && $\{I\}\cup \Cl_3^{\langle[\pi,2\pi]\rangle}$
  && $G\setminus\{-I\}$ && $G\setminus\{I\}$
  &\cr\h
& $\alpha+\beta+\gamma>2\pi$ &&
                $\big(\{-I\}\cup \Cl_3^{]\alpha+\beta+\gamma-2\pi,\pi]}\big)^\co$
  && $G\setminus\{-I\}$ && $G$
  &\cr\h
\noalign{\hrule}
}} }
\botcaption{ Table \tabProdABG } Triple and some quadruple products
of conjugacy classes in $SL_2(\R)$ involving $\Cl_3^\alpha\Cl_3^\beta\Cl_3^\gamma$
with $0<\alpha,\beta,\gamma<\pi$.
\endcaption
\endinsert

\midinsert
\centerline{\epsfxsize=110mm\epsfbox{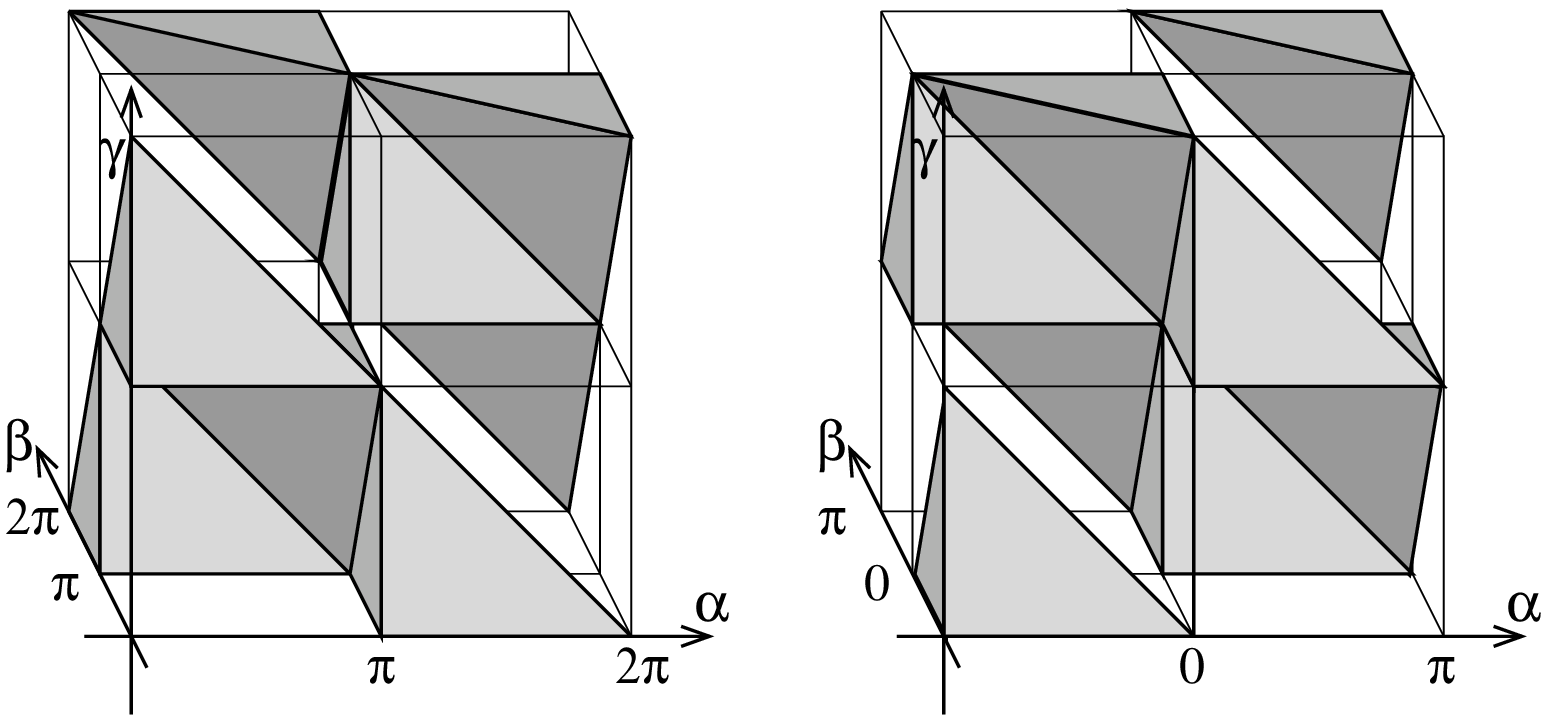}}
\botcaption{ Figure \figABG } The sets
$\big\{(\alpha,\beta,\gamma)\in\,(\,{]}0,2\pi{[}\,\setminus\{\pi\})^3 \mid
       I\in\Cl_3^\alpha\Cl_3^\beta\Cl_3^\gamma\big\}$
\break
 and
$\big\{(\alpha,\beta,\gamma)\in\,(\,{]}-\pi,\pi{[}\,\setminus\{0\})^3 \mid
       I\in\Cl_3^\alpha\Cl_3^\beta\Cl_3^\gamma\big\}$.
\endcaption
\endinsert

\head \sectTr. Range of the trace on the product of two classes
\endhead

Let $\Phi$ be as in (\eqIsom); see the introduction.

\proclaim{ Lemma \lemSU }
a). $SU(1,1)=\{\smat a & b \\ \bar b & \bar a \esmat\mid a,b\in \C,
a\bar a-b\bar b=1\}$.

\smallskip
b). $\Phi(\Cl_3^\alpha)$ is the conjugacy class of
$\left(\smallmatrix\lambda&0\\0&\bar\lambda\endsmallmatrix\right)$,
где $\lambda=e^{\alpha i}$.
\endproclaim

\demo{ Proof } Let $A=\smat a&b\\c&d\esmat$, $ad-bc=1$. Then
$A\in SU(1,1)$ if and only if $A^*J=JA^{-1}$. We have:
$$
   A^*J=\smat \bar a&\bar c\\ \bar b&\bar d\esmat\smat1&0\\0&-1\esmat=
   \smat\bar a&-\bar c\\ \bar b&-\bar d\esmat,
 \quad
   JA^{-1}=\smat1&0\\0&-1\esmat \smat d&-b\\-c&a\esmat
            = \smat d&-b\\ c&-a\esmat \qed
$$
\enddemo

\proclaim{ Lemma \lemEV }
Let $A\in\Cl_3^{[0,\pi]}$ and $B\in\overline{\Cl_4^+}$.
Then $AB\not\in\Cl_3^{[\pi,2\pi]}\cup\{I,-I\}$ unless \break
$A$ and $B$ both belong to $\Cl_2$ and have common eigenvector.
\endproclaim

\demo{Proof}
Let $v$ be a real eigenvector of $B$ that is not an eigenvector of $A$.
The condition $B\in\overline{\Cl_4^+}$ implies that the corresponding eigenvalue
$\lambda$ is positive. We have $v\wedge Av\ge 0$ by Proposition \propClass.
Moreover, $v\wedge Av\ne 0$ since $v$ is not an eigenvector of $A$.
Hence $v\wedge ABv = \lambda v\wedge Av>0$ and the result follows
from Proposition \propClass.
\qed\enddemo

\proclaim{ Lemma \lemAB} Let $0<\alpha,\beta<\pi$. Then:
\smallskip

(a). $\{\tr(AB)\,|\,A\in\Cl_3^{\alpha},\,B\in\Cl_3^{\beta}\,\}
   = \,{]}-\infty,\,2\cos(\alpha+\beta)]$;

\smallskip
(b). If $A\in\Cl_3^{\alpha}$, $B\in\Cl_3^{\beta}$, and
$\tr(AB)=2\cos(\alpha+\beta)$, then $AB\in\Cl_3^{\alpha+\beta}$.

\endproclaim

\demo{Proof} Let $A\in\Cl_3^{\alpha}$ and $B\in\Cl_3^{\beta}$.
We may assume by Lemma \lemSU\ that
$\Phi(A)=\left(\smallmatrix\lambda&0\\0&\bar\lambda\endsmallmatrix\right)$
and
$\Phi(B)=Q\left(\smallmatrix\mu&0\\0&\bar\mu\endsmallmatrix\right)Q^{-1}$ with
$\lambda=e^{\alpha i}$, $\mu=e^{\beta i}$, and
$Q=\left(\smallmatrix a & b\\ \bar b &\bar a\endsmallmatrix\right)$ with
$a\bar a-b\bar b=1$. Then we have
$$
\split
   \tr(AB) &= (\lambda\mu+\bar\lambda\bar\mu)a\bar a -
             (\lambda\bar\mu+\bar\lambda\mu)b\bar b
    = 2(1+b\bar b)\cos(\alpha+\beta) - 2b\bar b\cos(\alpha-\beta)
\\
    &= 2\cos(\alpha+\beta) - 4b\bar b\sin\alpha\sin\beta
\endsplit
$$
and the result easily follows.
\qed\enddemo

\proclaim{ Lemma \lemTwoThree }
Let $0<\alpha<\pi$. Then
$\{\tr(AB)\,|\,A\in\Cl_3^{\alpha},\,B\in\Cl_2^{++}\,\}
   = \,{]}-\infty,\,2\cos\alpha{[}$ and
$\{\tr(AB)\,|\,A\in\Cl_3^{\alpha},\,B\in\Cl_2^{+-}\,\}
   = \,{]}2\cos\alpha,\infty{[}$.
\endproclaim

\demo{Proof}  Let $A\in\Cl_3^{\alpha}$ and $B\in\Cl_2^{+\pm}$.
Let us fix a quadratic form invariant under $A$ and
choose a positively oriented orthonormal base $(e_1,e_2)$ such that
$e_2$ is an eigenvector of $B$. In this base, the matrices of the
corresponding operators are:
$A'=\left(\smallmatrix \cos\alpha &-\sin\alpha \\
                       \sin\alpha & \cos\alpha \endsmallmatrix\right)$ and
$B'=\left(\smallmatrix1&0\\\pm p&1\endsmallmatrix\right)$ with $p>0$,
thus $\tr A'B'=2\cos\alpha\mp p\sin\alpha$.
\qed\enddemo

\proclaim{ Lemma \lemTwoTwo }
 $\{\tr(AB)\,|\,A,B\in\Cl_2^{++}\} =
  \{\tr(AB)\,|\,A,B\in\Cl_2^{+-}\} = \,{]}-\infty,\,2]$ and
 $\{\tr(AB)\,|\,A\in\Cl_2^{++},\,B\in\Cl_2^{+-}\,\} = \,[2,\,\infty{[}$.
\endproclaim

\demo{Proof} We consider only the case $A,B\in\Cl_2^{++}$ (the other two cases
are similar). Let $A=\smat 1&0\\1&1\esmat$ and
$B=\smat a&b\\c&d\esmat\in\Cl_2^{++}$. Then $b\le0$ (see Table \tabClasses)
and $\tr AB = a+b+d = b+\tr B=b+2$. Moreover, $b$ can attain any non-positive value.
Indeed, consider the matrices $B_1=\smat 1&b\\0&1\esmat$ ($b<0$) and $B_2=A$.
\qed\enddemo

\proclaim{ Lemma \lemCFour } Let $A\in\Cl_4$. Then for any $t_1,t_2\in\R$ there
exist matrices $B,C\in G$ such that $\tr B=t_1$, $\tr C=t_2$, and $AB=C$.
\endproclaim

\demo{Proof}
Without loss of generality we may assume that $A=\diag(\lambda,\lambda^{-1})$,
$|\lambda|>1$. Let $B=\left(\smallmatrix a&b\\c&d\endsmallmatrix\right)$.
Then we have $\tr B=a+d$ and $\tr AB=\lambda a+\lambda^{-1}d$.
Thus, it is enough to find $a$ and $d$ from the simultaneous equations
$a+d=t_1$, $\lambda a+\lambda^{-1}d=t_2$ and then to find $b$ and $c$ such that
$bc = ad-1$.
\qed\enddemo


\head \sectDoubleProd. Double products of conjugacy classes
\endhead

\proclaim{ Lemma \lemABG }
Let  $0<\alpha,\beta<\pi$. Then $\Cl_3^\alpha \Cl_3^\beta$ is as shown
in Table \tabProd.
\endproclaim

\demo{Proof}
Let $A\in\Cl_3^\alpha$,
$B\in\Cl_3^\beta$, and let $C=AB$.

It follows from Lemma \lemAB\ that the range of $\tr C$ is
as required. So, it remains to show that the conjugacy class of $C$
is uniquely determined by the trace. This is evidently so when $C\in\Cl_4$.
Let us consider all the other cases.

Case 1. $\alpha+\beta<\pi$.
It is clear that $C\not\in\{\pm I\}$.

Case 1.1. $C\in\Cl_2^{\eps,\delta}$, $\eps,\delta=\pm1$.
We have $\eps=-1$ by Lemma \lemAB,
so it remains to show that $\delta\ne-1$. Suppose that
$\delta=-1$, i.e., $C\in\Cl_2^{--}$.
We have $B = (-A^{-1})(-C)$ with $-A^{-1}\in\Cl_3^{\pi-\alpha}$ and $-C\in\Cl_2^{++}$.
Hence Lemma \lemTwoThree\ implies $2\cos\beta=\tr B<\tr(-A^{-1})=2\cos(\pi-\alpha)$
which contradicts the assumption that $\alpha+\beta<\pi$.

Case 1.2. $C\in\Cl_3$.
Let $C\in\Cl_3^\gamma$, $\gamma\in{]}0,2\pi{[\,}\setminus\{\pi\}$.
Then, by Lemma \lemAB, we have
$\cos\gamma={1\over2}\tr C\le\cos(\alpha+\beta)$, hence
$$
    \gamma\ge\alpha+\beta.                              \eqno(\eqABone)
$$
Thus it suffices to show that $\gamma$ cannot be $>\pi$.

Suppose that $\gamma>\pi$.
Without loss of generality
we may assume that $\alpha\ge\beta$.
We have $(-A^{-1})(-C)=B$ with
$-A^{-1}\in\Cl_3^{\pi-\alpha}$, $-C\in\Cl_3^{\gamma-\pi}$.
The both angles $\pi-\alpha$ and $\gamma-\pi$ are in ${]}0,\pi{[}$. Thus,
by Lemma \lemAB\ applied to the matrices $-A^{-1}$, $-C$, and $B$,
we have
$\cos\beta={1\over2}\tr B\le\cos((\pi-\alpha)+(\gamma-\pi))
=\cos(\gamma-\alpha)$.
Combining this inequality with (\eqABone), we obtain $\gamma-\alpha\ge2\pi-\beta$.
Thus $\gamma\ge2\pi+\alpha-\beta$ which contradicts our assumptions
$\alpha\ge\beta$ and $\gamma<2\pi$.

\smallskip
Case 2. $\alpha+\beta=\pi$. Follows immediately from Lemma \lemAB.

\smallskip
Case 3. $\alpha+\beta>\pi$.
We have $C^{-1}=(-B^{-1})(-A^{-1})$ with
$-B^{-1}\in\Cl_3^{\pi-\beta}$,
$-A^{-1}\in\Cl_3^{\pi-\alpha}$, and $(\pi-\beta)+(\pi-\alpha)\in{]}0,\pi{[}$,
thus we reduce this case to Case 1.
\qed\enddemo

\proclaim{ Lemma \lemTwoThreeAB }
Let  $0<\alpha<\pi$. Then $\Cl_3^\alpha \Cl_2^{+\pm}$ is as shown
in Table \tabProd.
\endproclaim

\demo{Proof} Combine Lemma \lemTwoThree\ with Lemma \lemEV.
\qed\enddemo

\if01{ 
Due to Lemma \lemTwoThree\ it suffices to show that the trace
uniquely determines the conjugacy class
of any matrix $C=AB$ such that $A\in\Cl_3^\alpha$ and $B\in\Cl_2^{+\pm}$.
It is enough to consider two cases.

\smallskip
Case 1. $C\in\Cl_3$.
The result follows from Lemma \lemABG.

\smallskip
Case 2. $C\in\Cl_2$. We consider only the case $B\in\Cl_2^{+-}$
(the case $B\in\Cl_2^{++}$ is similar). Let $\vec x$ be the eigenvector of $B$.
Then $C\vec x=AB\vec x=A\vec x$, hence $\vec x\wedge C\vec x=\vec x\wedge A\vec x>0$
by Proposition \propClass(b). We have also $\tr C=2$ by Lemma \lemTwoThree.
Hence $C\in\Cl_2^{++}$ by Proposition \propClass(a).
}\fi 

\proclaim{ Lemma \lemTwoTwoAB }
$\Cl_2^{++}\Cl_2^{++}$,
$\Cl_2^{++}\Cl_2^{+-}$, and
$\Cl_2^{+-}\Cl_2^{+-}$ are as shown
in Table \tabProd.
\endproclaim

\demo{Proof}
$\Cl_2^{++}\Cl_2^{++}$ is as required due to Lemma \lemTwoTwo\ combined with
Lemma \lemEV. We obtain $\Cl_2^{+-}\Cl_2^{+-}$ from $\Cl_2^{++}\Cl_2^{++}$
by passing to the inverse matrices.
\if01{ 
Let us compute $X=\Cl_2^{++}\Cl_2^{++}$ (the computation of
$\Cl_2^{+-}\Cl_2^{+-}$ reduces to this case by passing to the inverse matrices).
By Lemma \lemTwoTwo, it is enough to show that: $\pm I\not\in X$ (which is evident),
$X\cap\Cl_3^{]\pi,2\pi[}=\varnothing$ (which follows from Lemma \lemTwoThreeAB),
and $X\cap\Cl_2^{*-}=\varnothing$. Let us prove the latter assertion.
Let $A,B\in\Cl_2^{++}$ and $AB\in\Cl_2$.
Let $\vec x$ be an eigenvector of $B$. If $\vec x$ is an eigenvector of $A$, then
$AB\in\Cl_2^{++}$. Otherwise we have $\vec x\wedge AB\vec x=\vec x\wedge A\vec x>0$
by Proposition \propClass(a) whence $AB\in\Cl_2^{*+}$, again by Proposition \propClass(a).
}\fi 
The computation of $\Cl_2^{++}\Cl_2^{+-}$ is immediate from Lemma \lemTwoTwo\
combined with the observation that $\smat 1&0\\1&1\esmat\smat 1&0\\-t&1\esmat$
belongs to $\Cl_2^{++}$, $\{I\}$, or $\Cl_2^{+-}$ according to the sign of $1-t$.
\qed\enddemo

\proclaim{ Lemma \lemCFourAB }
The product of $\Cl_4^\lambda$ with any other conjugacy class is as
in Table \tabProd.
\endproclaim

\demo{Proof}
Let $\lambda>1$, and $\Cl$ be any conjugacy class. Let $X$ be the set
that should coincide with $\Cl_4^\lambda\Cl$ according to Table \tabProd.
The fact that the product $\Cl_4^\lambda\Cl$ is contained in $X$,
is either evident or follows from Lemma \lemEV.
Let us prove the inverse inclusion. Let $A\in\Cl_4^\lambda$, $B_0\in\Cl$, and
$C_0\in X$. We assume that $B_0,C_0\ne\pm I$ (otherwise everything is evident).
We set $t_1=\tr B_0$ and $t_2=\tr C_0$. By Lemma \lemCFour, we can choose
matrices $B$ and $C$ such that
$\tr B=t_1$, $\tr C=t_2$, and $AB=C$. By passing to inverse matrices
if necessary, we may assume that $B\sim B_0$ (note that $A^{-1}\sim A$).
We have $\tr C=\tr C_0$ whence $C$ belongs either to the class of $C_0$ or to the class
of $C_0^{-1}$. However only one of these two
classes may be contained in $X$ (see Table \tabProd) which completes the proof.
\qed\enddemo

All double products of conjugacy classes are computed in Lemmas \lemABG--\lemCFourAB.
Using them, one easily computes triple and quadruple products as well; see
the subsequent sections.
%


\head \sectTripleProd. Triple products of conjugacy classes
\endhead

Due to previous computations we have:

\def\<{\langle}
\def\>{\rangle}
$$
\split
  \Cl_2^{++}\Cl_2^{++}\Cl_2^{++}&=
    \Cl_3^{[0,\pi]\rangle}\Cl_2^{++}
    =\big(\Cl_3^{]0,\pi[}\cup
     \Cl_2^{++}\cup
    -\Cl_2^{+-}\cup
      -\Cl_4^+\big)\Cl_2^{++}\\
  &=
    \Cl_3^{]0,\pi]\rangle} \cup
    \Cl_3^{[0,\pi]\rangle} \cup
    -\overline{\Cl_4^+} \cup
    -\Cl_3^{\langle[0,\pi]\rangle} = \{I\}^\co,
\endsplit
$$
$$
\split
  \Cl_2^{++}\Cl_2^{++}\Cl_2^{+-}&=
    \Cl_3^{[0,\pi]\>}\Cl_2^{+-}
    =\big(\Cl_3^{]0,\pi[} \cup
      \Cl_2^{++} \cup
     -\Cl_2^{+-} \cup
     -\Cl_4^+\big)\Cl_2^{+-}\\
  &=
      \Cl_3^{\<[0,\pi[}   \cup
       \overline{\Cl_4^+} \cup
      -\Cl_3^{\<[\pi,2\pi]} \cup
      -\Cl_3^{\<[\pi,2\pi]\rangle}
   =  \big(\{-I\}\cup \Cl_3^{[\pi,2\pi[}\big)^\co,
\endsplit
$$
$$
\split
   \Cl_3^\alpha\, \Cl_2^{++}\Cl_2^{++}&=
    \Cl_3^{]\alpha,\pi]\>}\Cl_2^{++}=
    \big(
       \Cl_3^{]\alpha,\pi[} \cup -\Cl_2^{+-} \cup -\Cl_4^+\big)\Cl_2^{++}
\\&
   = \Cl_3^{]\alpha,\pi]\>} \cup -\overline{\Cl_4^+} \cup -\Cl_3^{\<[0,\pi]\>}
   = \big(\{I\}\cup \Cl_3^{[0,\alpha]}\big)^\co,
\endsplit
$$
$$
\split
   \Cl_3^\alpha\, \Cl_2^{++}\Cl_2^{+-}&=
    \Cl_3^{]\alpha,\pi]\>}\Cl_2^{+-}=
    \big(
       \Cl_3^{]\alpha,\pi[} \cup -\Cl_2^{+-} \cup -\Cl_4^+\big)\Cl_2^{+-}
\\&
   = \Cl_3^{\<[0,\pi[} \cup -\Cl_3^{\<[\pi,2\pi]} \cup -\Cl_3^{\<[\pi,2\pi]\>}
   = \Cl_3^{\<[0,\pi]\>},
\endsplit
$$
$$
   \Cl_4^\lambda\, \Cl_2^{++}\Cl_2^{++}=
    \Cl_4^\lambda  \Cl_3^{[0,\pi]\>}
   =\Cl_4^\lambda \big(\Cl_3^{[0,\pi]}\cup -\Cl_4\big)
   =\Cl_3^{\<[0,\pi]\>}\cup -\{-I\}^\co = \{I\}^\co,
$$
$$
   \Cl_4^\lambda\, \Cl_2^{++}\Cl_2^{+-}
    =\Cl_4^\lambda\overline{\Cl_4^+}
     =\Cl_4^\lambda\big(\Cl_4^+ \cup \<\text{a subset of ${\Cl_4}^\co$}\>
     \big)
     = \{-I\}^\co,
$$
$$
  \Cl_4^\lambda\,\Cl_3^\alpha\,\Cl_2^{++}
 =\Cl_4^\lambda \Cl_3^{]\alpha,\pi]\>}
 =\Cl_4^\lambda\big(-\Cl_4^+ \cup \<\text{a subset of ${\Cl_4}^\co$}\>\big)
 =\{I\}^\co,
$$
$$
  \Cl_4^\lambda\,\Cl_3^\alpha\,\Cl_3^\beta
 =\Cl_4^\lambda\big(-\Cl_4^+ \cup \<\text{a subset of ${\Cl_4}^\co$}\>\big)
 =\{I\}^\co,
$$
$$
  \Cl_4^\lambda\, \Cl_4^\mu\,\Cl = G \qquad\qquad\text{for any non-scalar class $\Cl$}.
$$
If $\alpha+\beta<\pi$, then
$$
\split
  \Cl_3^\alpha\, \Cl_3^\beta\,\Cl_2^{++} &=
   \Cl_3^{[\alpha+\beta,\pi]\>}\Cl_2^{++} =
   \big( \Cl_3^{[\alpha+\beta,\pi[} \cup -\Cl_2^{+-}\cup -\Cl_4^+\big)\Cl_2^{++}
\\&=
   \Cl_3^{]\alpha+\beta,\pi]\>} \cup -\overline{\Cl_4^+}\cup -\Cl_3^{\<[0,\pi]\>}
   = \big(\{I\}\cup \Cl_3^{[0,\alpha+\beta]}\big)^\co.
\endsplit
$$
If $\alpha+\beta=\pi$, then
$$
  \Cl_3^\alpha\, \Cl_3^\beta\,\Cl_2^{++} =
   \big(\{-I\}\cup-\Cl_4^+\big)\Cl_2^{++} =
    -\Cl_2^{++}\cup -\Cl_3^{\<[0,\pi]\>} = \Cl_3^{\<[\pi,2\pi]\>}.
$$
If $\alpha+\beta>\pi$, then
$$
\split
  \Cl_3^\alpha\, \Cl_3^\beta\,\Cl_2^{++} &=
   \Cl_3^{\<[\pi,\alpha+\beta]}\Cl_2^{++} =
   \big( -\Cl_3^{]0,\alpha+\beta-\pi[} \cup -\Cl_2^{++}\cup -\Cl_4^+\big)\Cl_2^{++}
\\&=
   -\Cl_3^{]0,\pi]\>} \cup -\Cl_3^{[0,\pi]\>}\cup -\Cl_3^{\<[0,\pi]\>}
   = \Cl_3^{\<[\pi,2\pi]\>}.
\endsplit
$$
If $\alpha+\beta+\gamma<\pi$, then
$$
\split
  \Cl_3^\alpha\,\Cl_3^\beta\,\Cl_3^\gamma &=
   \Cl_3^{[\alpha+\beta,\pi]\>}\Cl_3^\gamma =
   \big( \Cl_3^{[\alpha+\beta,\pi-\gamma[}\cup \Cl_3^{\pi-\gamma}
   \cup \Cl_3^{]\pi-\gamma,\pi[}\cup -\Cl_2^{+-}\cup -\Cl_4^+\big)\Cl_3^\gamma
\\&=
  \Cl_3^{[\alpha+\beta+\gamma]\>}\cup
    (\{-I\}\cup-\Cl_4^+)\cup \Cl_3^{\<[\pi,\pi+\gamma[}\cup
   -\Cl_3^{\<[0,\gamma[} \cup -\Cl_3^{\<[0,\pi]\>}
\\&
   = \big(\{I\}\cup \Cl_3^{[0,\alpha+\beta+\gamma[}\big)^\co.
\endsplit
$$
If $\alpha+\beta+\gamma=\pi$, then
$$
\split
  \Cl_3^\alpha\,\Cl_3^\beta\,\Cl_3^\gamma &=
   \Cl_3^{[\alpha+\beta,\pi]\>}\Cl_3^\gamma =
   \big( \Cl_3^{\alpha+\beta}\cup \Cl_3^{]\alpha+\beta,\pi[}\cup
   -\Cl_2^{+-}\cup -\Cl_4^+\big)\Cl_3^\gamma
\\&=
   (\{-I\}\cup \Cl_4^-)\cup \Cl_3^{\<[\pi,\pi+\gamma[}\cup
    -\Cl_3^{\<[0,\gamma[} \cup -\Cl_3^{\<[0,\pi]\>}
   =\{-I\}\cup \Cl_3^{\<[\pi,2\pi]\>}.
\endsplit
$$
If $\pi<\alpha+\beta+\gamma<2\pi$ and $\alpha+\beta<\pi$, then
$$
\split
  \Cl_3^\alpha\,\Cl_3^\beta\,\Cl_3^\gamma &=
   \Cl_3^{[\alpha+\beta,\pi]\>}\Cl_3^\gamma =
   \big( \Cl_3^{[\alpha+\beta,\pi[}\cup -\Cl_2^{+-}\cup -\Cl_4^+\big)\Cl_3^\gamma
\\&=
   \Cl_3^{\<[\pi,\pi+\gamma[}\cup -\Cl_3^{\<[0,\gamma[} \cup -\Cl_3^{\<[0,\pi]\>}
    = \Cl_3^{\<[\pi,2\pi]\>}.
\endsplit
$$
If $\pi<\alpha+\beta+\gamma<2\pi$ and $\alpha+\beta=\pi$, then
$$
  \Cl_3^\alpha\,\Cl_3^\beta\,\Cl_3^\gamma =
   \big(\{-I\}\cup-\Cl_4^+\big)\Cl_3^\gamma
    = -\Cl_3^\gamma\cup - \Cl_3^{\<[0,\pi]\>}= \Cl_3^{\<[\pi,2\pi]\>}.
$$

All other triple products reduce to these by
passing to the inverse matrices (and maybe changing the sign).
For example, if $\pi<\alpha+\beta+\gamma<2\pi$ and $\alpha+\beta>\pi$, then
$$
   \pi<(\pi-\alpha)+(\pi-\beta)+(\pi-\gamma)<2\pi,\qquad
    (\pi-\alpha)+(\pi-\beta)<\pi,
$$
hence
$$
  \Cl_3^\alpha\,\Cl_3^\beta\,\Cl_3^\gamma =
  (-\Cl_3^{\pi-\alpha})^{-1}(-\Cl_3^{\pi-\beta})^{-1}(-\Cl_3^{\pi-\gamma})^{-1}
   = \big(-\Cl_3^{\<[\pi,2\pi]\>}\big)^{-1}=\Cl_3^{\<[\pi,2\pi]\>}.
$$


\head \sectQuadrProd. Quadruple products. End of proof of Theorem \thMain
\endhead

We see in Table \tabProd\ that any triple product of non-scalar conjugacy classes contains
$\Cl_4$, and the product of $\Cl_4$ with any other non-scalar class
contains $\{I,-I\}^\co$. Thus, to compute the product $X$ of four non-scalar classes,
it is enough to check if $I$ and $-I$ belongs to it.
In its turn, to decide whether or not $\pm I\in X$, it is enough to check if
the inverse of one of these four classes multiplied by $\pm1$
belongs to the product of the others.

For example, let us check that $X:=\Cl_3^\alpha\Cl_3^\beta\Cl_3^\gamma\Cl_2^{++}=\{-I\}^\co$
when $\alpha+\beta+\gamma=\pi$ (see Table \tabProdABG).
Indeed,
$$
\xalignat2
    \big(\Cl_2^{++}\big)^{-1}=\Cl_2^{+-}&\in\Cl_3^\alpha\Cl_3^\beta\Cl_3^\gamma
     =\{-I\}\cup\Cl_3^{\<[\pi,2\pi]\>}, &&\text{hence $I\in X$},
\\
    -\big(\Cl_2^{++}\big)^{-1}=\Cl_2^{-+}&\not\in\Cl_3^\alpha\Cl_3^\beta\Cl_3^\gamma
     =\{-I\}\cup\Cl_3^{\<[\pi,2\pi]\>}, &&\text{hence $-I\not\in X$}.
\endxalignat
$$

As another example, let us compute the product
$\Cl_3^\alpha\Cl_3^\beta\Cl_3^\gamma\Cl_3^\delta$ which we denote by $X$.
We consider only the case when $\alpha+\beta+\gamma+\delta\le 2\pi$ because
the other case, when this sum is in the range $[2\pi,4\pi[$, is
reduced to this one by passing to the inverse matrices.
Let $Y=\Cl_3^\alpha\Cl_3^\beta\Cl_3^\gamma$. We have
$$
    Y=\cases
         \Big(\{I\}\cup\Cl_3^{[0,\alpha+\beta+\gamma[}\Big)^\co, &\alpha+\beta+\gamma<\pi,\\
          \{-I\}\cup\Cl_3^{\<[\pi,2\pi]\>},    &\alpha+\beta+\gamma=\pi,\\
           \Cl_3^{\<[\pi,2\pi]\>},    &\alpha+\beta+\gamma>\pi.
      \endcases
$$
Therefore
$\big(\Cl_3^\delta\big)^{-1}=-\Cl_3^{\pi-\delta}\in Y$ whence $I\in X$, and since
$-\big(\Cl_3^\delta\big)^{-1}=\Cl_3^{\pi-\delta}$, we have
$$
   -I\in X \;\Leftrightarrow\; \Cl_3^{\pi-\delta}\in Y
           \;\Leftrightarrow\; \pi-\delta\not\in{]}0,\alpha+\beta+\gamma{[}
           \;\Leftrightarrow\; \pi-\delta\ge\alpha+\beta+\gamma
$$
whence $X=G$ if $\alpha+\beta+\gamma+\delta\le\pi$, and $X=\{-I\}^\co$
if $\pi<\alpha+\beta+\gamma+\delta\le2\pi$.


\Refs

\ref\no\refAW
\by     S.~Agnihotri, C.~Woodward
\paper  Eigenvalues of products of unitary matrices and quantum Schubert
        calculus \jour Math. Research Letters \vol 5 \yr 1998 \pages 817--836
\endref

\ref\no\refAH
\by    Z.~Arad, M.~Herzog (eds.)
\book  Products of conjugacy classes in groups \bookinfo Lecture Notes in Math. 1112
\publ  Springer-Verlag  \publaddr Berlin, Heidelberg, N.Y., Tokyo \yr 1985
\endref

\ref\no\refACC\by E.~Artal Bartolo, J.~Carmona Ruber, J.~I.~Cogolludo Agust\'\i n
\paper Effective invariants of braid monodromy 
\jour Trans. Am. Math. Soc. \vol 359 \issue 1 \pages 165--183 \yr 2007 \endref

\ref\no\refBelkale
\by     P.~Belkale
\paper  Local systems on $\Bbb P^1-S$ for $S$ a finite set
        \jour Compos. Math. \vol 129 \yr 2001 \pages 67--86
\endref

\ref\no\refB\by J.~I.~Brenner
\paper Covering theorems for nonabelian simple groups, IV
\jour J\~n\B an\B abha Ser. A\vol 3 \yr 1973 \pages 77--84 \endref

\ref\no\refFW\by E.~Falbel, R.~A.~Wentworth
\paper On products of isometries of hyperbolic space
\jour Topology Appl. \vol 156 \yr 2009 \pages 2257--2263 \endref

\ref\no\refGP\by O.~Garc\'\i a-Prada, M.~Logares, V.~Mu\~noz
\paper Moduli spaces of parabolic $U(p,q)$-Higgs bundles 
\jour Q. J. Math. \vol 60 \issue 2\pages 183--233 \yr 2009 \endref 

\ref\no\refGordeev
\by     N.L.~Gordeev
\paper  Products of conjugacy classes in perfect linear groups. Extended covering number
\jour   Zapiski Nauchn. Semin. POMI \vol 321 \yr 2005 \pages 67--89
\transl Reprinted in:
\jour   J. Math. Sci. \vol 136 \yr 2006 \pages 1867--3879
\endref


\ref\no\refKarni
\by     S.~Karni
\paper  Covering number of groups of small order and sporadic groups
\inbook Ch.~3 in [\refAH] \pages 52--196
\endref

\ref\no\refLSh
\by     M.W.~Liebeck, A.~Shalev
\paper  Diameter of finite simple groups: sharp bounds and applications
\jour   Ann. of Math. \vol 154 \yr 2001 \pages 383--406
\endref

\ref\no\refJKTR\by S.~Yu.~Orevkov
\paper Quasipositivity test via unitary representations of braid groups and its
applications to real algebraic curves 
\jour J. Knot Theory Ramifications \vol 10 \issue 7 \pages 1005--1023 \yr 2001 \endref

\ref\no\refAFST\by S.~Yu.~Orevkov
\paper Products of conjugacy classes in finite unitary groups $GU(3,q^2)$ and $SU(3,q^2)$
\jour Ann. Fac. Sci. de Toulouse. Math\'ematiques (6)
\vol 22 \issue 2 \yr 2013 \pages 219--251 \endref

\ref\no\refP\by J.~Paupert
\paper Elliptic triangle groups in $PU(2,1)$, Lagrangian triples and momentum maps
\jour Topology \vol 46 \yr 2007 \pages 155--183 \endref

\ref\no\refPW\by J.~Paupert, P.~Will
\paper Involution and commutator length for complex hyperbolic isometries
\jour Michigan Math. J. \vol 66 \yr 2017 \pages 699--744 \endref

\ref\no\refS\by C.~Simpson
\paper Products of matrices \inbook
Differential geometry, global analysis, and topology (Halifax, NS, 1990)
\pages 157--185
\bookinfo CMS Conf. Proc., 12
\publ Amer. Math. Soc. \publaddr Providence, RI \yr 1991.\endref

\endRefs
\enddocument